\documentclass[12pt]{article}
\usepackage{hyperref}
\usepackage{graphicx}
\usepackage{amssymb,amsmath}
\newtheorem{Th}{Theorem}
\newtheorem{Def}{Definition.}
\newenvironment{Proof}
{\par\noindent{\bf Proof.}}
{\hfill$\scriptstyle\square$}
\newtheorem{remark}{Remark.}

\begin{document}
	
	\title{LOCAL LIMIT THEOREMS FOR MARKOV CHAINS WITH TREND COMPONENT OF LINEAR GROWTH \footnote{This research was supported by the Scientific Fund of  NRU HSE}}
	\date{}
	\author{V. Konakov \\ \href{mailto:VKonakov@hse.ru}{VKonakov@hse.ru} \and A. Markova \\ \href{mailto:anna.r.markova@gmail.com}{anna.r.markova@gmail.com} \\ \and
		Laboratory of Stochastic Analysis and its Applications, \\
		National research university Higher School of Economics\\
		Shabolovka, 26, Moscow, 119049 Russia }
	\maketitle
	
	\emph{\bf {Abstract}}: We consider a sequence of Markov chains weakly convergent to a diffusion. We suppose that a drift term contains a linearly increasing component. The usual parametrix method fails because of this unbounded drift term. We show how to modify the parametrix method to obtain local limit theorems for this case.  
	\\
	

\section{Introduction}

The aim of the introduction is to explain the origin of the problem considered in this research. It is necessary to mention that the introduction is not intended to be a full and detailed review.

The importance of researching the discretization of stochastic differential equations (SDE) results from the fact that every known approximate modeling method for SDE is based upon some discretization. It is common knowledge that the solution to SDE in closed form can be obtained only for a very limited number of stochastic differential equations subclasses, consequently different discretization methods (the Euler-Maruyama scheme, the Milstein scheme; the higher order stochastic Taylor expansions schemes) draw attention of a large number of researchers. Herewith an evident issue emerges: whether discretization schemes reproduce the asymptotic properties of the initial process (such as   mixing rate, large deviation principle, etc.). Another important issue is the following: how “close” is the designed approximation to the solution of the initial SDE in this or that way. Our research concerns the second issue. Historically the first considered scheme was the sequence of Markov chains defined on a grid with the time step tending to zero; weak convergence of measures (transition probabilities) to the transition probability of a limited diffusion process was being studied. The first general results concerning weak convergence for such a scheme were obtained by A.V. Skorokhod in 1961 [1]. Skorokhod$'$s results deal with quite a general class of Markov chains and processes which could have jumps. The most general weak convergence results for continuous diffusion were obtained in the monograph by D. Stroock and S. Varadhan [2]. They developed an approach based upon the solution of the so called “martingale problem”. We should emphasize once again that these first results concerned weak convergence of measures. The modern theory of weak convergence of probability measures has been developed for many decades and is associated primarily with the names of A.N. Kolmogorov, J.L. Doob, M. Donscker, Yu.V. Prokhorov, A.V. Skorokhod, L. Le Cam ans S. Varadhan.

Let us now assume that transition probabilities of Markov chains and of limiting diffusion process are absolutely continuous w.r.t. the Lebesgue measure. In this case we naturally ask a question of when the convergence of these transition densities occurs, in other words when the corresponding local limit theorem holds true. In order to answer this question V. Konakov and S. Molchanov [3] introduced a discrete version of the parametrix method which was later on specified and generalized in the series of publications by V. Konakov and E. Mammen [4-7]. It is worth mentioning that the parametrix method has been known for a long time in the theory of differential equations. It was introduced by E. Levi in 1907 [8, 9] and then developed in the studies by A. Friedman [10], A. Il$'$yin, A. Kalashnikov and O. Oleinik [11] and many others. 
But this version of the parametrix method did not fit our aims. H.P. McKean and I.M. Singer [12] introduced a modification of the parametrix method in 1967 which, as it turned out, allows for a discrete version and makes it possible to develop a new method of obtaining local limit theorems for transition densities of Markov chains sequence, which are weakly convergent to the limiting diffusion process. One of the essential conditions for proving these results was the condition of bounded coefficients of drift and diffusion. Limitation was necessary for the convergence of designed parametrix series. This condition narrowed the sphere of the applicability of the acquired results and prevented from considering a number of important specific models. The goal of this research is the description of the procedure which allows excluding linearly increasing trend component and reducing the problem to the already investigated one with bounded drift and diffusion coefficients. This procedure is applied to both diffusion and Markov chains. For diffusion similar procedure of trend exclusion has been applied  in the paper by F. Delarue and S. Menozzi [13] in the context of obtaining two-sided estimates of transition density for some degenerate  Kolmogorov type SDE$'$s. In reference to Markov chains this procedure, so far as we are aware, is new. The essence of the procedure is simple: the trend growth should be compensated by a backtrack along the trajectories of an ODE system, which is acquired by means of deleting Brownian component of SDE. The stochastic differential and the corresponding SDE of this compensated process can be easily written by the Ito$'$s formula. This SDE has a bounded drift coefficient. A similar procedure is performed for Markov chains, where instead of a differential equation a difference equation is used and the backtrack occurs not along the trajectory of a differential equation, but along its Euler$'$s broken line. After that, having applied the known results to a bounded case, it is possible to return back to the initial problem by a simple transformation and get local limit theorems for the initial model. Later on  the authors are also intending to consider a more general case of increasing trend with bounded gradient and the case of  unbounded diffusion coefficient

\section{Necessary information from the differential equations theory and the difference equations theory}

We consider the linear homogeneous vector differential equation with variable coefficients:
\begin{equation} \label{1}
	x^{'}=A(t)x, x \in R^{d}, A(t): R^{d} \mapsto R^{d}  ,t \in [0,1],
\end{equation}
where $A(t)$ is continuous matrix valued function.

We consider the set $X$ of all solutions to the equation (\ref{1}) defined on the interval $[0,1]$. The set $X$ is a linear vector space consisting of functions $phi:[0,1] \mapsto R^{d}$

Along with the vector differential equation (\ref{1}) we consider the matrix differential equation:
\begin{equation} \label{2}
X^{'}=A(t)X        
\end{equation}

We say that a matrix function $\Phi(t)$ is a solution of (\ref{2}) on the interval $[0,1]$ if $\Phi(t)$ is continuously differentiable for  $t \in [0,1]$ и $\Phi^{'}(t)=A(t)\Phi(t), t \in [0,1]$.

The following theorem establishes the connection between the solutions of equations (\ref{1}) and (\ref{2}).

\begin{Th}\cite[Theorem 2.20, p. 33]{14} %
Let $A(t)$ be a continuous  $n \times {n}$ be a continuous  $[0,T]$ and let  $\Phi(t)$ be an  $n \times n$ matrix-valued function with the columns  $ \phi_{1}(t), \phi_{2}(t),…,\phi_{n}(t)$:
$$
\Phi(t)=[\phi_{1}(t), \phi_{2}(t),…,\phi_{n}(t)], t \in [0,T].
$$

Then $\Phi$ is the solution of the matrix differential equation (\ref{2}) on $[0,T]$ if and only if each column $\phi_{i}$ is a solution of the vector differential equation (\ref{1}) on $[0,T]$, $i=1,2,...,n$. Moreover, if $\Phi$ is the solution of the matrix equation (\ref{2}), then
$$
x(t)=\Phi(t)c
$$
is the solution of the vector differential equation (\ref{1}) for any  $n \times 1$ vector of constants $c$.
\end{Th}

We introduce the following definitions.

\begin{Def}
	. The fundamental system of solutions of the differential equation (\ref{1}) is a basis of the vector space $X$.
\end{Def}
\begin{Def}
	A matrix whose columns form a fundamental system of solutions is called the fundamental matrix of the differential equation  (\ref{1}).
\end{Def}

We consider the linear homogeneous difference equation of order  $k$:
\begin{equation} \label{3}
	x(s+k)+a_{1}(s)x(s+k-1)+...+a_{k}(s)x(s)=0.
\end{equation}
The set of solutions of the difference equation (\ref{3}) ) is also a vector space. By analogy with the definition for the differential equation -  \emph{the fundamental matrix of the difference equation} 
is a matrix whose columns form a
\emph{fundamental system of solutions of the difference equation (\ref{3})}.And the fundamental system of solutions of the difference equation forms a basis in the vector space of all solutions of the difference equation.

We should note that it is convenient to take the identity matrix of the appropriate dimension as the fundamental matrix of the differential or difference equation at the initial time $t=0$. In this case $x(t)=\Phi(t)c$ is the solution of the vector equation (\ref{1}) with the initial conditions $x(0)=\Phi(0)c=c$.

\section{The procedure of excluding of the linear trend component for the diffusion model and for the Markov chain}

We consider the following diffusion model:
\begin{equation} \label{4}
dY=\{b(t)Y+m(t,Y)\}dt+\sigma(t,Y)dB(t), Y(0)=x \in R^{d}, t \in [0,1],                
\end{equation}
where $B(t)$ is the standard Brownian motion. The interval $[0,1]$ has been chosen for convenience and can be replaced by any bounded interval.

We also consider a triangle array of Markov chains with the same initial conditions as in the model (\ref{4}):
\begin{equation} \label{5}
	\begin{array}{lcl}
X_{n} \left( \frac{k+1}{n} \right) =X_{n} \left( \frac{k}{n} \right) + \frac{1}{n} \left\lbrace b_{n}\left(\frac{k}{n}\right)X_{n}\left(\frac{k}{n}\right)+m_{n} \left(\frac{k}{n},X_{n}\left(\frac{k}{n}\right)\right)\right\rbrace+\frac{1}{\sqrt{n}}\varepsilon_{n}\left(\frac{k+1}{n}\right), \\
\\
X_{n}(0)= x \in R^{d}, k=0,1,2,...,n.                                                	
	\end{array}
\end{equation}
We make the standard Markov assumption for innovations  $\varepsilon_{n}$, namely, the random variable  $\varepsilon_{n}\left(\frac{k+1}{n}\right)$ given the past $X_n\left(\frac{i}{n}\right)=x(i),i=0,1,2,..k$ depends only on the value of the process $x(k)$ at the last moment of time $\frac{k}{n}$ and has a conditional density $q_{n,\frac{k}{n},x(k)}(\cdot)$ taken from a family of densities $q_{n,t,x}(\cdot)$ dependent depending on a triplet $(n,t,x) \in N \times [0,1] \times R^{d}$. As concerns the family of densities $q_{n,t,x}(\cdot)$ and the coefficients of the equation (\ref{4}) the following assumptions are supposed to be satisfying:

\begin{enumerate}
	\item
	$a(t,x)=\sigma(t,x)\sigma^{T}(t,x)$ is a symmetric positively definite matrix such that $c\le \theta^{T}a(t,x)\theta \le C$ для всех $\theta$ for any $\mid \theta \mid=1, x \in R^{d}, t \in [0,1]$.
	
	\item
	The matrix functions $b_n(t)$, $b(t)$ are continuous on $[0,1]$. The functions $a(t,x)$ and $m(t,x)$ and their first derivatives w.r.t. t and x are continuous and bounded uniformly in $(t,x), x \in R^{d}, t \in [0,1]$ and Lipchitz continuous in the space variable $x$ with Lipchitz constant independent on $t$. Moreover, the second derivatives  $\frac{\partial^{2}a(t,x)}{\partial x_i \partial x_j}, 1 \le i, j \le d$ exist and satisfy H\"{o}lder condition in $x$ with the constant independent on $t$. 
	
	\item
	$\int q_{n,t,x}(z)zdz=0$, $\int q_{n,t,x}(z)zz^{T}dz \stackrel{\triangle}{=}a_{n}(t,x)$, $x \in R^{d}, t \in [0,1]$.
	
	\item
	There exist a positive integer  $S^{'}$ and a function $\psi(x):R^{d} \mapsto R$, $\sup_{x \in R^{d}} \mid \psi (x) \mid \ < \infty$ such that $\int \mid x \mid ^{s} \psi (x)dx < \infty$, where $S=2dS^{'}+4$ and for every  n large enough and $n$, $z \in R^{d}$:
	\begin{eqnarray*}
		\begin{array}{lcl}
			\mid D_{z}^{v} q_{n,t,x}(x) \mid \le \psi (z) , \mid v \mid =0,1,2,3,4, \\
			 \mid D_{x}^{v}	q_{n,t,x}(x) \mid \le \psi (z) , \mid v \mid =0,1,2,
			 	
		\end{array}
	\end{eqnarray*}
where $x \in R^{d}, t \in [0,1]$.
		
	\item
For a sequence $\bigtriangleup_{n} \mapsto 0$:
\begin{eqnarray*}
	\begin{array}{lcl}
\sup_{x \in R^{d}, t \in [0,1]} \mid m_{n} (t,x)-m(t,x) \mid = O \left(\bigtriangleup_{n}\right), \\
\sup_{x \in R^{d}, t \in [0,1]} \mid a_{n} (t,x)-a(t,x) \mid = O \left(\bigtriangleup_{n}\right), \\
\sup_{t \in [0,1]} \mid b_{n} (t,x)-b(t,x) \mid = O \left(\bigtriangleup_{n}\right).
	\end{array}
\end{eqnarray*}
\end{enumerate}

The existence of transition density for the model (\ref{5}) follows immediately from the model assumptions. The existence of transition density for the diffusion model (\ref{4}) may be derived from the H\"{o}rmander theory \cite{15} but under stronger conditions on the coefficients. The existence of transition density in this model can also be proved by the parametrix method under weaker conditions for coefficients than the conditions required for the H\"{o}rmander theory. We consider the diffusion model  (\ref{4}) first. If the model$'$s coefficients are bounded, then the existence of the transition density follows from the paper by Il$'$yin, Kalashnikov and Oleinik \cite{11} or from the paper by Konakov and Mammen \cite{4}, where the parametrix was constructed for this equation. But in the model (\ref{4}) the trend is unbounded and linearly increasing. To apply the parametrix method to the this model, we propose the following procedure: exclude the linear component and consider a new model with bounded coefficients; then apply the parametrix method to this new model with bounded coefficients and go back to the original model with linearly increasing trend.

We consider the linear system of ordinary differential equations (ODE):
\begin{equation} \label{6}
y^{'}(t)=b(t)y(t), y(0)=x \in R^{d}, t \in [0,1]                                             
\end{equation}

Let $\Phi(t)$ be the fundamental matrix corresponding to this system. We remind that for the system (\ref{6}) it means that all solutions may be prolongated to the whole interval $[0,1]$ and the fundamental matrix is the matrix whose columns are independent solutions of this system corresponding to the given initial conditions. We take the identity matrix: $\Phi(0)=I$, as the initial fundamental matrix. Clearly, the fundamental matrix is the solution of the equation $\Phi^{'}(t)=b(t)\Phi(t)$ and with our initial condition is nondegenerate on the whole interval $[0,1]$. The inverse matrix $\Phi^{-1}(t)$ satisfies the equation $[\Phi^{-1}(t)]^{'}=-\Phi^{-1}(t)b(t)$ and the initial condition $\Phi^{-1}(0)=I$. To exclude the linear component of the trend let us consider the process $\tilde{Y}(t)=f\left(t,Y(t)\right)$, where $f(t,y)=\Phi^{-1} (t)y$. The function $f(t,y)$ is continuous on $[0,1] \times R^{d}$ and has continuous partial derivatives $\frac{\partial f}{\partial t}$, $\frac{\partial f}{\partial y_{i}}$, so we can use the Ito$'$s formula to obtain the stochastic differential for the process $\tilde{Y}(t)$:
\begin{eqnarray*}
\begin{array}{lcl}
d\tilde{Y}(t)=d\left[ \Phi^{-1}(t)Y(t)\right]=\Phi^{-1}(t)dY(t)+d\Phi^{-1} (t)Y(t)= \\
=\Phi^{-1}(t) \left( \{ b(t)Y(t)+m \left( t,Y(t) \right) \} dt + \sigma \left( t,Y(t) \right) dB(t) \right) + [\Phi^{-1}(t)]^{'}dt \; Y(t)= \\
= \Phi^{-1}(t)b(t)Y(t)dt+\Phi^{-1}(t) m \left( t,Y(t) \right) dt+ \Phi^{-1}(t)\sigma \left( t,Y(t) \right) dB(t) - \\
- \Phi^{-1}(t) b(t) Y(t) dt = \Phi^{-1}(t)  m \left( t,Y(t) \right) dt + \Phi^{-1}(t) \sigma \left( t,Y(t) \right)  dB(t)= \\
=\tilde{m} \left( t, \tilde{Y}(t) \right) dt + \tilde{\sigma}  \left( t,\tilde{Y}(t) \right) dB(t),
\end{array}
\end{eqnarray*}
where $\tilde{m} \left( t, \tilde{Y}(t) \right)= \Phi^{-1}(t)  m \left( t, \Phi (t) \tilde{Y}(t) \right)$,  $\tilde{\sigma} \left( t, \tilde{Y}(t) \right)= \Phi^{-1}(t)  \sigma \left( t, \Phi (t) \tilde{Y}(t) \right)$.

We see that the process $\tilde{Y}(t)$ is the diffusion process satisfying SDE with the bounded drift $\tilde{m} \left( t, \tilde{Y}(t) \right)$ and with positively definite diffusion matrix $\tilde{\sigma} \left( t, \tilde{Y}(t) \right)$:
$$
d\tilde{Y}(t)=\tilde{m} \left( t, \tilde{Y}(t) \right) dt+\tilde{\sigma} \left( t, \tilde{Y}(t) \right) dB(t).
$$

Indeed:
\begin{eqnarray*}
\begin{array}{lcl}
\tilde{a}(t,y)=\tilde{\sigma} \left( t, \tilde{Y}(t) \right) [\tilde{\sigma} \left( t, \tilde{Y}(t) \right)]^{T}=\Phi^{-1}(t) \; \sigma \left( t, \Phi(t) y \right) [\Phi^{-1}(t) \; \sigma \left( t, \Phi(t) y \right)]^{T}= \\
= \Phi^{-1}(t) \sigma \left( t, \Phi(t) y \right) [\sigma \left( t, \Phi(t) y \right)]^{T} [\Phi^{-1}(t)]^{T} = \Phi^{-1}(t) \; a \left( t, \Phi(t) y \right) [\Phi^{-1}(t)]^{T}
\end{array}
\end{eqnarray*}
and, hence,	
$$
\theta ^{T} \tilde{\sigma} (t,y) [\tilde{\sigma} (t,y)]^ {T} \theta = 0 \Leftrightarrow \vartheta ^{T} \sigma (t,y) [\sigma (t,y)]^ {T} \vartheta=0, \vartheta = [\Phi^{-1}(t)]^ {T} \vartheta.
$$

It remains to use positive definiteness of the matrix $a(t,\Phi(t)y)=\sigma(t,\Phi(t)y) \times [\sigma(t,\Phi(t)y)]^{T}$.

The existence of the transition density $\rho_{\tilde{Y}} (t)$ of the process $\tilde{Y}(t)$ is proved using the parametrix method in the article by Konakov and Mammen \cite{4}. By the well known transformation formulas the transition density of the process $Y(t)=\Phi (t) \tilde{Y}(t)$ is equal to:
\begin{equation} \label{7}
\rho_{Y} (s,t,x,y)=\det [\Phi^{-1}(t)]\rho_{\tilde{Y}}(s,t,\Phi^{-1}(s)x,\Phi^{-1}(t)y)
\end{equation}

We consider the trend exclusion procedure for the model (\ref{5})), which is a discrete analogue of the procedure described above for diffusion equation. We consider the following difference equation without trend:
$$
\frac{x_{n} \left( (k+1)h \right) - x_{n}(kh)}{h} = b_{n}(kh) X_{n}(kh), \; x_{n}(0)=x
$$
on the grid $\Gamma=\{0,h,2h,…,nh=1\}, h=\frac{1}{n}$.

In the matrix notation:
$$
x_{n} \left( (k+1)h \right)=\left( I+h \; b_{n}(kh)  \right) x_{n}(kh), \; x_{n}(0)=x.
$$

Iterating we obtain:
\begin{eqnarray*}
	\begin{array}{lcl}
	x_{n}(h)=\left( I+h \; b_{n}(0) \right)x, \\
	x_{n}(2h)=\left( I+h \; b_{n}(h) \right)x_{n}(h)=\left( I+h \; b_{n}(h) \right) \left( I+h \; b_{n}(0) \right)x, \\
	... \\
	x_{n}(kh)=\Phi_{n}(kh)x,
	\end{array}
\end{eqnarray*}
where $\Phi_{n}(kh)=\left( I+h \; b_{n}\left( (k-1)h\right)\right) \Phi_{n} \left( (k-1)h \right) \Phi_{n} (kh)$ is the fundamental matrix in the theory of difference equations \cite{16}, a discrete analogue of the fundamental matrix $\Phi(t)$ is defined on the grid  $\{0,h,2h,…,nh=1\}$, with the initial condition $\Phi_{n}(0)=I$.

We define a new Markov chain:
$$
\tilde{X}_{n}(kh)=\Phi_{n}^{-1}(kh) \; X_{n}(kh), \tilde{X}_{n}(0)=x.
$$

By the well known formulas of density transformation, the transition density of the Markov chain  $X_{n}(kh)= \Phi_{n}(kh) \; \tilde{X}_{n} (kh)$  is equal to:   
\begin{equation} \label{8}
 	\rho_{X_{n}} (ih,jh,x,y)=\det [\Phi^{-1}(jh)]\rho_{\tilde{X}_{n}}(ih,jh,\Phi^{-1}(ih)x,\Phi^{-1}(jh)y)
\end{equation}

Then following the model (\ref{5}) we have:
\begin{eqnarray*}
\begin{array}{lcl}
\tilde{X}_{n} \left((k+1)h \right)= \Phi_{n}^{-1} \left( (k+1)h \right) \; X_{n} \left( (k+1)h \right) = \\
= \Phi_{n}^{-1} (kh) \left( I+ h \; b_{n}(kh) \right)^{-1} \{ \left( I+ h \; b_{n}(kh) \right)X_{n}(kh) +\\
+ h \; m_{n} \left( kh, \; X_{n}(kh) \right) +\sqrt{h} \varepsilon_{n}  \left(   (k+1)h \right)   \} = \\
=\Phi_{n}^{-1} (kh) X_{n}(kh) + h \; \Phi_{n}^{-1} (kh)\; \left(    I+ h \; b_{n}(kh) \right)^{-1} \; m_{n} \left(kh, X_{n}(kh) \right) + \\
+\sqrt{h} \; \Phi_{n}^{-1} (kh) \; \left( I+h \; b_{n} (kh) \right)^{-1} \; \varepsilon_{n}  \left(   (k+1)h \right)  = \\
=\tilde{X}_{n}(kh)+h \; \tilde{m}_{n} \left( kh, \; \tilde{X}_{n}(kh) \right) + \sqrt{h} \; \tilde{\varepsilon}_{n} \left( (k+1)h \right),
\end{array}
\end{eqnarray*}
where
\begin{eqnarray*}
\begin{array}{lcl}
	\tilde{m}_{n} \left( kh, \; \tilde{X}_{n} (kh)  \right) = \Phi_{n}^{-1} (kh) \; \left( I+h \; b_{n} (kh) \right)^{-1} \; m_{n} \left( kh, \; \Phi_{n} (kh) \; \tilde{X}_{n} (kh)  \right), \\
	\tilde{\varepsilon}_{n} \left( (k+1)h \right) = \Phi_{n}^{-1} (kh) \; \left( I+h \; b_{n} (kh) \right)^{-1} \; \varepsilon_{n} \left( (k+1)h \right).
\end{array}
\end{eqnarray*}

It follows from the definition of $\tilde{\varepsilon}_{n}$ that the random variable $\tilde{\varepsilon}_{n} \left(   (k+1)h \right)$ given the “past” $\tilde{X}_{n}(ih)=\tilde{x}(i), \; i=0,1,2,...,k$ depends only on the value of the process $\tilde{x}(k)$ at the last moment of time $kh$ and has the conditional density:
\begin{equation} \label{9}
\tilde{q}_{n,kh,\tilde{x}(k)}(z)=\det \; \Phi_{n} \left(   (k+1)h\right) \; q_{n,kh,\Phi_{n}(kh)\tilde{x}(k)} \left(  \Phi_{n} \left(  (k+1)h\right)z \right)
\end{equation}

This density is taken from the family of densities  $\det \; \Phi_{n} \left(   ([tn]+1)h\right) \times \tilde{q} _{n,t,\Phi_{n}\left(   ([tn]+1)h\right){x}(k)} \left(  \Phi_{n} \left(  ([tn]+1)h\right)z \right)$, which is dependent on the triplet $(n,t,x)$: $(n,t,x) \in N \times [0,1] \times R^{d}$. The densities $\tilde{q}_{n,kh,\tilde{x}(k)}(z)$ in (\ref{9}) satisfy conditions 3 and 4, formulated at the beginning of the paper, with allowance for the fact that $\phi (x)$ from condition 4 is replaced by $C \; \phi (x)$, where $C$ is a constant. Making the change of variables $v=\Phi \left( \frac{[tn]+1}{n} \right)z$, with $t=kh$, we have:
\begin{equation} \label{10}
	\begin{array}{lcl}
		\int \tilde{q}_{n,t,\tilde{x}}(z)dz=\det \Phi_{n} \left( \frac{[tn]+1}{n}  \right) \; \int q_{n,t,\Phi_{n} \left( \frac{[tn]}{n} \right) \tilde{x}}  \left( \Phi_{n} \left(  \frac{[tn]+1}{n} \right)  z \right)z dz = \\
		= \det \Phi_{n} \left( \frac{[tn]+1}{n}  \right) \; \det \Phi_{n}^{-1} \left( \frac{[tn]+1}{n}  \right) \;  \int q_{n,t,\Phi_{n} \left( \frac{[tn]}{n} \right) \tilde{x}}  (v)  \Phi_{n}^{-1} \left(  \frac{[tn]+1}{n} \right)  v dv = \\
		= \Phi_{n}^{-1} \left(  \frac{[tn]+1}{n} \right) \; \int q_{n,t,\Phi_{n} \left( \frac{[tn]}{n} \right) \tilde{x}}  (v)  v dv  = 0, \\
		\int  \tilde{q}_{n,t,\tilde{x}}(z) z_{i} z_{j} dz = \det \Phi_{n} \left( \frac{[tn]+1}{n}  \right) \; \int q_{n,t,\Phi_{n} \left( \frac{[tn]}{n} \right) \tilde{x}}  \left( \Phi_{n} \left(  \frac{[tn]+1}{n} \right)  z \right) z_{i} z_{j} dz = \\
		=  \int q_{n,t,\Phi_{n} \left( \frac{[tn]}{n} \right) \tilde{x}}  (v)  \left[ \Phi_{n}^{-1} \left(  \frac{[tn]+1}{n} \right)  v \right]_{i}  \left[ \Phi_{n}^{-1} \left(  \frac{[tn]+1}{n} \right)  v \right]_{j} dv = \\
		= \left\{ \Phi_{n}^{-1} \left(  \frac{[tn]+1}{n} \right) \; \int q_{n,t,\Phi_{n} \left( \frac{[tn]}{n} \right) \tilde{x}}  (v)  v v^{T} dv  \left[ \Phi_{n}^{-1} \left(  \frac{[tn]+1}{n} \right)    \right]^{T}  \right\}_{ij} = \\
		= \left\{ \Phi_{n}^{-1} \left(  \frac{[tn]+1}{n} \right) \; a_{n} \left(t, \Phi_{n} \left( \frac{[tn]}{n}  \right) \tilde{x}  \right) \left[ \Phi_{n}^{-1} \left(  \frac{[tn]+1}{n} \right)    \right]^{T}  \right\}_{ij} \stackrel{\triangle}{=} \tilde{a}_{n} (t, \tilde{x}) 
	\end{array}
\end{equation}

The vector function $\Phi_{n} (t) x$ coincides in the points $t=kh, k=0,1,2,...,n$ with the Euler broken line for the equation $y^{'} (t)=b_{n} (t)y(t),   y(0)=x \in R^{d}$. Therefore in the case of a diffusion process we compensate the increasing trend with the backtrack along the trajectories of a differential equation $y^{'} (t)=b_{n} (t)y(t),   y(0)=x \in R^{d}$, and in the case of a Markov chain with the backtrack along the Euler broken line of this equation. According to the well-known properties of the Euler broken lines \cite{17}, $\Phi_{n} \left( \left[ \frac{t}{h}  \right] \right) x \rightarrow \Phi (t)$ uniformly on the interval  $[0,1]$ and, given the properties described above, from (\ref{6}) and (\ref{10}) we have:
\begin{eqnarray*}	
	\tilde{a}_{n} (t, \tilde{x}) = \int \tilde{q}_{n,t,x} (z) z z^{T} \rightarrow \tilde{a} (t, x) = \tilde{\sigma} (t,x) [\tilde{\sigma} (t,x)]^{T}, \;\;\; n \rightarrow \infty.
\end{eqnarray*}
	
Let us now illustrate how the assertions obtained for the models with bounded trends are transformed into the corresponding assertions for the models with linear component in trends. For simplicity we consider a case of the density family $q_{n,t,x} (\cdot)$, independent on the parameter  $n$, in other words $q_{n,t,x} (\cdot) = q_{t,x} (\cdot)$. Let us assume that the conditions 1 – 4 listed above are satisfied for the family $q_{t,x} (\cdot)$ and the coefficients of the equation (\ref{4}) (condition 5 in our case is trivially fulfilled, as $m_{n} (t,x) \equiv m(t,x)$ и   $a_{n} (t,x) \equiv a(t,x)$). Then for $\tilde{m}(t,x) = \Phi^{-1} (t) ;\ m(t, \Phi (t)x)$ and $\tilde{\sigma}(t,x) = \Phi^{-1} (t) \; \sigma (t, \Phi (t)x)$  the conditions of Theorem 1.1 from V. Konakov and E. Mammen \cite{4} are fulfilled,  and so we have the estimate
\begin{equation} \label{11}
\sup_{x,y \in R^{d}} \left( 1+ ||y-x||^{2(S^{'}-1)} \right) ||p_{\tilde{X}_{n}} (0,1,x,y) - p_{\tilde{Y}} (0,1,x,y) ||= O (\frac{1}{\sqrt{n}})
\end{equation}

Using (\ref{7}) and (\ref{8}), we formulate the following result, which follows from (\ref{11}).

\begin{Th}
Let all the conditions 1)-4) listed above be satisfied. Then
\begin{eqnarray*}
		\sup_{x,y \in R^{d}} \left( 1+ ||\Phi
		 (1)y-x||^{2(S^{'}-1)} \right) ||\det \Phi (1) \; p_{X_{n}} (0,1,x,\Phi (1)y) - \\
		 -\det \Phi (1) \; p_{Y} (0,1,x,\Phi (1)y) ||= O (\frac{1}{\sqrt{n}}).
\end{eqnarray*}
If $b(t) \equiv b$, then
\begin{eqnarray*}
		\sup_{x,y \in R^{d}} \left( 1+ ||\Phi ^{-1}
		(1)y-x||^{2(S^{'}-1)} \right) ||p_{X_{n}} (0,1,x,y) - \\
		-p_{Y} (0,1,x,y) ||= O (\frac{1}{\sqrt{n}}).
\end{eqnarray*}
\end{Th}

\begin{Proof}
	Let us apply (\ref{11}) to the points $x$ and $\Phi ^{-1} (1)y$,  then we get
\begin{equation} \label{12}
		\begin{array}{lcl}
			\sup_{x,y \in R^{d}} \left( 1+ ||\Phi ^{-1}
			(1)y-x||^{2(S^{'}-1)} \right) ||p_{\tilde{X}_{n}} (0,1,x,\Phi ^{-1} (1)y) - \\
			-p_{\tilde{Y}} (0,1,x,\Phi ^{-1} (1)y) ||= O (\frac{1}{\sqrt{n}}).
		\end{array}
\end{equation} 
From (\ref{7}) and (\ref{8}) we have
\begin{equation} \label{13}
	 p_{\tilde{Y}} (s,t,x,z) = \det \Phi (t) p_{Y} (s,t, \Phi (s)x, \Phi (t)z),
\end{equation} 
\begin{equation} \label{14}
  	 p_{\tilde{X}_{n}} (ih,jh,x,v) = \det \Phi_{n} (jh) p_{X_{n}} (ih,jh, \Phi _{n} (ih)x, \Phi _{n} (jh)v),
\end{equation} 

Substituting in (\ref{13}) and (14)  $s=0,t=1,i=0,j=n,x$ and $z= \Phi ^{-1} (1) y, \; v=\Phi_{n}^{-1} (y)$,  we get from (\ref{12}):
\begin{equation} \label{15}
	\begin{array}{lcl}
		\sup_{x,y \in R^{d}} \left( 1+ ||\Phi
		(1)y-x||^{2(S^{'}-1)} \right) ||\det \Phi (1) \; p_{X_{n}} (0,1,x,\Phi (1)y) - \\
		-\det \Phi (1) \; p_{Y} (0,1,x,\Phi (1)y) ||= O (\frac{1}{\sqrt{n}}).
\end{array}
\end{equation} 

Let $b(t) \equiv b$ .  Then $\Phi (1) =e^{b}$ is a matrix exponent, $\Phi_{n} (1)= \left(  I+ \frac{b}{n} \right)^{n}$ and for $n$ large enough
\begin{equation} \label{16}
|| \Phi (1) - \Phi _{n} (1)  ||\le e^{a} - \left( 1+ \frac{a}{n} \right)^{n} \le \frac{a^{2} e^{a}}{n},
\end{equation} 
where  $a=||b||$. The first inequality in  (\ref{16}) The first inequality in \cite{17}, p. 98, to prove the second inequality it is enough to find the sign of the derivative near the point $x=0$ of the function $f(x) \stackrel{\triangle}{=} a^{2} e^{a} x- e^{a}+ (1+ax)^{1/x} $. Besides from Lemma 3.1 and 3.2 in \cite{4} we have an estimate
$$
p_{\tilde{Y}} (s,t,x,y) \le C e^{-C||y-x ||^{2}},
$$
so for the transition density $p_{Y} (0,1,x,y)$ we get
\begin{equation} \label{17}
p_{Y} (0,1,x,y) = \det \Phi ^{-1} p_{\tilde{Y}} (0,1,x, \Phi ^{-1} (1) y) \le C e^{-C ||\Phi^{-1} (1) y - x||^{2}}.
\end{equation} 
The second assertion of the theorem now follows from (\ref{15}) - (\ref{17}).

\end{Proof}

\begin{remark}
	As	$\; C_{1} ||y- \Phi (1) x ||\le ||\Phi^{-1} (1) y - x ||\le C_{2} ||y- \Phi (1) x||,$
the assertions of the theorem can also be written with the factor $\left( 1+||y - \Phi (1) x||^{2(S^{'}-1)} \right)$ instead of $\left( 1+||\Phi ^{-1} (1) y -  x||^{2(S^{'}-1)} \right)$. So the nonuniform estimate of the convergence rate in this theorem is achieved either by pulling $y$ back, or by pushing $x$ forward.
\end{remark}

\section{The example of excluding linear component of the trend of the diffusion model}

Consider the model presented in B. Koo, O. Linton \cite{18}:
\begin{equation} \label{18}
 d X_{t}=\{ \beta (t) (a(t)-X_{t}) \} dt + \sigma (t, X_{t}) dB_{t}, \; \; \; \; X_{0}=x \in R^{d},
\end{equation} 
where $B_{t}, t \geq 0$ is the standard Brownian motion.

Consider a system of Linear Ordinary Differential Equations (LODE): $x^{'} (t) = -\beta (t) x(t), \; y(0)=x \in R^{d}$ and its fundamental matrix $\Phi (t): \Phi^{'} (t) = - \beta (t) \Phi (t), \; \Phi (0)=I$, where $I$ is the identity matrix. Consider the process $\tilde{X}_{t} = \Phi^{-1} (t) X_{t}$ and by the Ito$'$s Lemma \cite{19} we obtain the stochastic differential:
\begin{eqnarray*}
 d \tilde{X}_{t}=d \left[ \Phi^{-1} (t) X_{t}  \right] = \Phi^{-1} (t) d X_{t} + d \Phi^{-1} (t) X_{t} = 
 \Phi^{-1} (t)\beta (t) a(t) dt + \\
 + \Phi^{-1} (t) \sigma (t, X_{t}) dB_{t}
 =\tilde{m} (t, \tilde{X}_{t}) dt + \tilde{\sigma} (t, \tilde{X}_{t}) dB_{t}
\end{eqnarray*}
where $\tilde{m} (t, \tilde{X}_{t}) = \Phi^{-1} (t)\beta (t) a(t), \; \tilde{\sigma} (t, \tilde{X}_{t}) = \Phi^{-1} (t) \sigma (t, \Phi (t) X_{t})$.

So process $\tilde{X}_{t}$ can be represented as:
\begin{eqnarray*}
		\tilde{X}_{t}=\tilde{X}_{0} + \int_{0}^{t} \tilde{m} (s, \tilde{X}_s) ds + \int_{0}^{t} \tilde{\sigma} (s, \tilde{X}_s) d B_{s} = x+ \int_{0}^{t} \Phi ^{-1} (s) \beta (s) a(s) d s + \\+
		\int_{0}^{t} \Phi ^{-1} (s) \sigma (s, \Phi (s) \tilde{X}_{s}) d B_{s}. 
\end{eqnarray*}

Consider the one-dimensional case with constant coefficients: $a(t)\equiv a, \; \beta (t) \equiv \beta, \; \sigma (t) \equiv \sigma$. Model (\ref{18}) then corresponds to the model of the evolution of the interest rate by Vasicek \cite{20}:
\begin{equation} \label{19}
	d X_{t}=\{ \alpha \beta- \beta X_{t} \} dt + \sigma dB_{t}, \; \; \; \; X_{0}=x \in R
\end{equation} 

For the model (\ref{19}) there is explicit solution of the SDE. So the procedure of excluding trend of the diffusion equation is not particularly interesting, however, allows us to trace the correctness of the method.

The solution of the SDE (\ref{19}) can be represented in \cite{20}:
\begin{equation} \label{20}
	X_{t}=x \; \exp (-\beta t) +a(1-\exp (-\beta t)) + \sigma \; \exp (-\beta t) \int_{0}^{t} \exp (\beta s) dB_{s}.
\end{equation} 

Consider getting the solution to SDE (\ref{19}) using the procedure of excluding the linear component of the trend.

The fundamental matrix for SDE (\ref{19}) has the form: $\Phi (t) =\exp (-\beta t) $. Then process $\tilde{X}_{t}=\Phi ^{-1} (t) X_{t} = \exp (\beta t) X_{t}$ can be represented as:
$$
\tilde{X}_{t} = x + a (\exp (\beta t) -1) + \sigma \int_{0}^{t} \exp (\beta s) d B_{s}.
$$

The solution obtained by the exclusion of the linear component of the trend, takes the form:
\begin{eqnarray*}
X_{t} = \Phi (t) \left( x+ a(\exp (\beta t) -1) +\sigma \int_{0}^{t} \exp (\beta s) dB_{s}  \right) = x \; \exp (-\beta t) + \\
+ a(1- \exp (-\beta t)) + \sigma \; \exp (-\beta t) \int_{0}^{t} \exp (\beta s) dB_{s},
\end{eqnarray*}
that is coincides with formula (\ref{20}).

The method of excluding of the linear component of the trend is also applicable to more difficult models, such as a modified model of the evolution of the interest rate by Cox-Ingersoll-Ross \cite{21}, and to its extension - the Hull-White model \cite{22}. Besides the models of interest rates, the trend exclusion method can be applied to the stochastic volatility model \cite{23} in the case when the function of volatility contained in the profitability equation is bounded. The Cox-Ingersoll-Ross model is differs from the above Vasicek model in function of volatility $\sigma (t, X_{t})\equiv \sigma \sqrt{X_{t}}$, procedure of excluding the linear component of the trend for it is analogous to the above procedure for the Vasicek model.

Consider the two-dimensional case for the Heston stochastic volatility model \cite{23}:
\begin{align} \label{21}
\left\{  \begin{array}{cc} d S_{t} = \mu S_{t} dt + f(v_{t}, S_{t}) dB_{t}^{1}, \\ d v_{t} = k(\theta - v_{t}) dt + \xi g(v_{t}) dB_{t}^{2}, \end{array} \right.
\end{align}
where $B_{t}=(B_{t}^{1},B_{t}^{2})$  is the standard Brownian motion, $f(v_{t},S_{t}), g(v_{t})$ are bounded functions.

We represent a system of SDE \cite{21} in matrix form
\begin{align*}
\left(	\begin{array}{cc} d S_{t} \\ d v_{t} \end{array}	\right)
&=
\left[  \left(	\begin{array}{cc} \mu & 0 \\ 0 & -k \end{array}	\right)
\left(	\begin{array}{cc} S_{t} \\ v_{t} \end{array}	\right) 
+
\left(	\begin{array}{cc} 0 \\ k \theta \end{array}	\right) \right] dt
+
\left(	\begin{array}{cc} f (v_{t}, S_{t}) & 0 \\ 0 & \xi g(v_{t}) \end{array}	\right)
\left(	\begin{array}{cc} dB_{t}^{1} \\ dB_{t}^{2} \end{array}	\right).
\end{align*}

Consider a system of LODE:
$$
x^{'} (t) =  \left(	\begin{array}{cc} \mu & 0 \\ 0 & -k \end{array}	\right) x(t).
$$

The fundamental matrix for this system has the form:
$$
\Phi (t) =  \left(	\begin{array}{cc} \exp (\mu t) & 0 \\ 0 & \exp(-kt) \end{array}	\right).
$$

Then the process
$$
\left(	\begin{array}{cc} \tilde{S}_{t} \\ \tilde{v}_{t} \end{array}	\right)
=
\left(	\begin{array}{cc} \exp (-\mu t) & 0 \\ 0 & \exp(kt) \end{array}	\right)
\left(	\begin{array}{cc} S_{t} \\ v_{t} \end{array}	\right).
$$
can be represented as:
\begin{eqnarray*}
\left(	\begin{array}{lc} \tilde{S}_{t} \\ \tilde{v}_{t} \end{array}	\right)
=
\left(	\begin{array}{cc} S_{0} \\ v_{0} \end{array}	\right)
+
\left(	\begin{array}{cc} 0 \\ \theta(\exp(kt)-1) \end{array}	\right)
+
\int_{0}^{t}
\left(	\begin{array}{cc} \exp (-\mu s) & 0 \\ 0 & \exp(ks) \end{array}	\right)
\times \\ \times
\left(	\begin{array}{cc} f (v_{s}, S_{s}) & 0 \\ 0 & \xi g(v_{s}) \end{array}	\right) d B_{s}.
\end{eqnarray*}

Then the solution of the SDE (\ref{21}), obtained by exclusion of the linear trend component takes the form:
\begin{eqnarray*}
	\begin{array}{lcl}
	\left(	\begin{array}{lc} S_{t} \\ v_{t} \end{array}	\right)
	=
	\left(	\begin{array}{cc} \exp (\mu t) & 0 \\ 0 & \exp(-kt) \end{array}	\right) 
	\left(	\begin{array}{lc} S_{0} \\ v_{0} \end{array}	\right)
	+
	\left(	\begin{array}{cc} 0 \\ \theta(1- \exp(-kt)) \end{array}	\right)
	+ \\ +
	\left(	\begin{array}{cc} \exp (\mu t) & 0 \\ 0 & \exp(-kt) \end{array}	\right) 
	\int_{0}^{t}
	\left(	\begin{array}{cc} \exp (-\mu s) & 0 \\ 0 & \exp(ks) \end{array}	\right)
	\left(	\begin{array}{cc} f (v_{s}, S_{s}) & 0 \\ 0 & \xi g(v_{s}) \end{array}	\right) d B_{s}.
	\end{array}
\end{eqnarray*}

\section{Conclusion}
It is known \cite{11}, \cite{24} that the parametrix method and its discrete analogue \cite{4}, \cite{5} assume the boundedness of the drift and diffusion coefficients. Meanwhile quite a number of important models have non-limited drift coefficient, in particular the models corresponding to stochastic  recurrent  estimation procedures have a linearly increasing drift coefficient.

More precisely, Markov chains and limiting  diffusion processes with linear trend component occur in recurrent estimation procedures based upon the Robbins-Monroe method. A series of results concerning weak convergence of recurrent estimation procedures to finite-dimensional distributions of some limiting  diffusion process was proved in monograph \cite{25} (Theorem 6.3, Chapter 6; Theorems 3.1 and 5.2, Chapter 8). These theorems assume the existence of densities, that is why the natural question is: whether not only weak convergence, but also density convergence holds, that is whether the corresponding local limit theorem is true? The parametrix method combined with the trend exclusion method makes a positive answer possible. This application of the method will be the subject of a separate publication.

This research was aimed at introducing a procedure which would allow for excluding the linearly increasing trend component and reduce the problem to the previously studied one with bounded trend. The suggested method is applicable to more general models with trends having bounded gradient, but the formulas are less intuitive and Euler broken lines are drawn locally, so in this case we can speak of a local limit theorem in a short time. This will also be the topic of our next research.


\begin{thebibliography}{10}

\bibitem{1}
{\it A. V. Skorokhod,}
Studies in the Theory Random Processes // Dover Pubns, 1982.
\bibitem{2}
{\it D. W. Stroock, S. R. S. Varadhan,} Multidimensional diffusion processes // Grundlehren der Mathematischen Wissenschaften 233. Springer (Berlin), 1979.
\bibitem{3}
{\it V. Konakov, S. Molchanov,} On the convergence of Markov chains to diffusion processes // Theory Probab. and Math. Stat., v.31, Kiev Univ., 1984.
\bibitem{4}
{\it V. Konakov, E. Mammen,} Local limit theorems for transition densities of Markov chains converging to diffusions // Prob. Th. Rel. Fields, 117:551–587, 2000.
\bibitem{5}
{\it V. Konakov, E. Mammen,} Local approximations of Markov random walks by diffusions // Stoch. Proc. Appl., 96(1), pp. 73-98, 2001.
\bibitem{6}
{\it V. Konakov, E. Mammen,} Edgeworth type expansions for transition densities of Markov chains converging to diffusions // Bernoulli: a journal of mathematical statistics and probability. 2005. Vol. 11. No. 4. P. 591-641.
\bibitem{7}
{\it V. Konakov, E. Mammen,}  Small time Edgeworth-type expansions for weakly convergent nonhomogeneous Markov chains // Probability Theory and Related Fields. 2009. Vol. 143. No. 1. P. 137-176.
\bibitem{8}
{\it	E. E. Levy,} Sulle equazioni lineari alle derivate parziali totalmente ellittiche // Rendiconti della Reale Accademia dei Lincei, Classe di Scienze Fisiche, Matematiche, Naturali, Serie V, 16 (12): 932–938, 1907.
\bibitem{9}
{\it	E. E. Levy,} Sulle equazioni lineari totalmente ellittiche alle derivate parzial // Rendiconti del Circolo Matematico di Palermo, 24(1): 275-317, 1907.
\bibitem{10}
{\it	A. Friedmsn,} Partial differential equations of parabolic type // Prentice-Hall, 1964.
\bibitem{11}
{\it	A. M. Ilyin, A. S. Kalashnikov, and O. A. Oleynik,} Linear second-order partial differential equations
of the parabolic type // Journal of Mathematical Sciences, Vol. 108, No. 4, 2002.
\bibitem{12}
{\it	H. P. McKean, I. M. Singer,} Curvature and the eigenvalues of the Laplacian // J. Differential Geometry 1, 43–69, 1967.
\bibitem{13}
{\it	F. Delarue, S. Menozzi,} Density Estimation for a Random Noise Propagating through a Chain of Differential Equations // J. Func. Anal., 259, 2010, №6, 1577-1630.
\bibitem{14}
{\it	W. Kelley, A. Peterson,} The Theory of Differential Equations Classical and Qualitative // Prentice Hall, 2004.
\bibitem{15}
{\it	D. Nualart,} The Malliavin Calculus and Related Topics // Springer, 2006.
\bibitem{16}
{\it	S. Elaydi,} An Introduction to Difference Equations // Springer Science+Business Media, Inc, 2005.
\bibitem{17}
{\it  V.I. Arnold,} Ordinary Differential Equations // The MIT Press, 1978.
\bibitem{18}
{\it	B. Koo, O. Linton,} Semiparametric Estimation of Locally Stationary Diffusion Models // Discussion paper No. EM/2010/551, August 2010.
\bibitem{19}
{\it B. Oksendal,} Stochastic differential equations. An Introduction with Applications // Springer, 2000.
\bibitem{20}
{\it	O. Vasicek,} An Equilibrium Characterization of the Term Structure // Journal of Financial Economics 5 (2), pp. 177–188, 1977.
\bibitem{21}
{\it	J.C. Cox, J.E. Ingersoll, S.A Ross,} A Theory of the Term Structure of Interest Rates // Econometrica, 53, pp.:385-407, 1985.
\bibitem{22}
{\it	J. Hull, A. White,} Pricing interest-rate derivative securities // The Review of Financial Studies, Vol 3, No. 4, pp. 573–592, 1990.
\bibitem{23}
{\it	S. L. Heston,} A Closed-Form Solution for Options with Stochastic Volatility with Applications to Bond and Currency Options // The Review of Financial Studies,  Volume 6, number 2, pp. 327—343, 1993.
\bibitem{24}
{\it	V. D. Konakov,} Parametrix method for diffusions and Markov chains // Preprint, Moscow State University, 2012 (in Russian).
\bibitem{25}
{\it	M. Nevel'son, R. Khas'minskii,} Stochastic approximation and recurrent estimation // M, Nauka, 1972 (in Russian).
	
\end{thebibliography}
\end{document}